\documentclass[12pt,reqno]{amsart}
\usepackage{amscd,graphicx}
\usepackage{amssymb}
\usepackage{color}
\usepackage{tikz-cd}
\usepackage{mathtools}
\newcommand\SmallMatrix[1]{{%
  \tiny\arraycolsep=0.3\arraycolsep\ensuremath{\begin{bmatrix}#1\end{bmatrix}}}}

\newcommand*{\matminus}{%
  \leavevmode
  \hphantom{0}%
  \llap{%
    \settowidth{\dimen0 }{$0$}%
    \resizebox{1.1\dimen0 }{\height}{$-$}%
  }%
}

\usepackage[colorlinks,
    linkcolor={red!50!black},
    citecolor={blue!50!black},
    urlcolor={blue!80!black}]{hyperref}
\newcommand{\wh}{\widehat}

\newcommand{\sm}{\smallsetminus}

\newcommand{\al}{\alpha}
\newcommand{\be}{\beta}

\newcommand{\ep}{\varepsilon}
\renewcommand{\th}{\theta}

\renewcommand{\phi}{\varphi}
\newcommand{\si}{\sigma}

\newcommand{\Ga}{\Gamma}

\newcommand{\De}{\Delta}

\newcommand{\cO}{\mathcal O}

\newcommand{\Hom}{\operatorname{Hom}}

\newcommand{\conj}{\operatorname{\!conj}}

\newtheorem{theorem}{Theorem}

\theoremstyle{remark}

\begin{document}

\title[The ${SU(2)}$ Casson-Lin invariant of the Hopf link]
{The $\boldsymbol{SU(2)}$ Casson-Lin invariant of the Hopf link}
\author{Hans U. Boden}
\address{Mathematics \& Statistics, McMaster University, Hamilton, Ontario} 
\email{boden@mcmaster.ca}
\thanks{The first author was supported by a grant from the Natural Sciences and Engineering Research Council of Canada.}

\author{Christopher M. Herald}
\address{Mathematics \& Statistics, University of Nevada Reno, Reno Nevada}
\email{herald@unr.edu}


\date{\today}
\begin{abstract}
We compute the $SU(2)$ Casson-Lin invariant for the Hopf link and determine the sign in the formula of Harper and Saveliev relating this invariant to the linking number.
\end{abstract}
\maketitle

The Casson--Lin   knot  invariant   was defined by Lin \cite{Lin}, and then extended to  an invariant $h_2(L)$ of   2-component links by Harper and Saveliev \cite{HS1}, who showed 
that $h_2(L) = \pm lk (\ell_1,\ell_2)$, the linking number of  $L=\ell_1 \cup \ell_2$, up to an overall sign.
The purpose of this note is to determine the sign in that formula, establishing the following.

\begin{theorem} \label{thm-main}
If $L=\ell_1 \cup \ell_2$ is an oriented 2-component link in $S^3$, then its Casson-Lin invariant 
$h_2(L) = -lk (\ell_1,\ell_2)$.
\end{theorem}

Recall from \cite{HS1} that the Casson-Lin invariant $h_2(L)$ is  a signed count of certain projective $SU(2)$ representations of the link group.  This   paper is devoted to proving Theorem \ref{thm-main}. 
The proof of Proposition 5.7 in \cite{HS1} shows that the sign in the above formula is independent of $L$,  and thus Theorem \ref{thm-main} will follow from a single computation. 

To that end, we will determine the Casson-Lin invariant for the right-handed  Hopf link.
Since there is  just one irreducible projective $SU(2)$ representation of the link group, up to conjugation, it suffices to determine    the sign associated to this   one   point.  

We identify $$SU(2) = \{x+yi+zj+wk \mid |x|^2+|y|^2+|z|^2+|w|^2=1\}$$ with the group of unit quaternions and consider the conjugacy class $$C_i = \{ yi+zj+wk \mid |y|^2+|z|^2+|w|^2=1\} \subset SU(2)$$ of purely imaginary unit quaternions. Notice that $C_i$ is diffeomorphic to $S^2$ and  coincides with the set of  $SU(2)$ matrices of trace zero.

 Let $L$ be an oriented link in $S^3$, represented as the closure of an $n$-strand  braid $\si \in B_n$. We follow Conventions 1.13 on p.17 of \cite{KT} for writing geometric braids $\si$ as words in the standard generators $\si_1,\ldots, \si_{n-1}$.  In particular, braids are oriented from top to bottom and $\si_i$ denotes a right-handed crossing in which the $(i+1)$-st strand crosses over the $i$-th strand.  The braid group $B_n$ gives a faithful right action on the free group $F_n$ on $n$ generators, and here we follow the conventions in \cite{BHa} for associating an automorphism of $F_n$ to a given braid $\si \in B_n$, which we write as $x_i \mapsto x_i^\si$ for $i=1,\ldots, n$. Note that each braid $\si \in B_n$ fixes the product $x_1 \cdots x_n$.  

A standard application of the Seifert-Van Kampen theorem shows that the link complement $S^3 \sm L$ has fundamental group
$$\pi_1(S^3 \sm L) = \langle x_1, \dots, x_n \mid  x_i^\si = x_i, \, i=1, \dots, n  \rangle.$$  
We can therefore identify representations in $\Hom(\pi_1 (S^3 \sm L), SU(2))$ with fixed points in $\Hom(F_n, SU(2))$ under the induced action of the braid $\si$. We further identify $\Hom(F_n, SU(2))$ with $SU(2)^n$ by associating $\varrho$ with the $n$-tuple $(X_1,\ldots, X_n) = (\varrho(x_1), \ldots, \varrho_1(x_n))$.  Note that $\si \colon SU(2)^n \to SU(2)^n$ is equivariant with respect to conjugation, so that fixed points come in whole orbits. 

If $\ep=(\ep_1,\ldots, \ep_n)$ is an $n$-tuple with $\ep_i = \pm 1$ and $\ep_1 \cdots \ep_n = 1,$ then projective $SU(2)$ representations can be identified with fixed points in $\Hom(F_n, SU(2))$ under the induced action of $\ep \si$, which also preserves
the product $X_1 \cdots X_n$   and is conjugation equivariant.   The Casson-Lin invariant $h_2(L)$ is then defined as a signed count of orbits of 
fixed points of $\ep \si$ for a suitably chosen $n$-tuple $\ep=(\ep_1,\ldots, \ep_n)$.   The choice is made so that the resulting projective representations $\al$ all have $w_2(\al) \neq 0$ and as such do not lift to $SU(2)$ representations, and it has the consequence that for all fixed points $\al$ of $\ep \si$, each $\al(x_i)$ is a traceless  $SU(2)$ element.

 We therefore restrict our attention to the subset of traceless representations, which are elements $\varrho \in \Hom(F_n, SU(2))$ with $\varrho(x_j) \in C_i$ for $j=1,\ldots, n.$  Define $f\colon C_i^n \times C_i^n \to SU(2)$ by setting $$f(X_1, \dots, X_n, Y_1, \dots, Y_n)=(X_1 \cdots X_n) (Y_1 \cdots Y_n)^{-1}.$$ 
 We obtain an orientation on $f^{-1}(1)$ by applying the base-fiber rule, using the product orientation on $C_i^n \times C_i^n$ and the standard orientation on the codomain of $f$.  The quotient  $f^{-1}(1)/\conj $
 is then oriented by another application of the base-fiber rule, using the standard orientation on $SU(2)$. This step uses the fact that, if $\ep=(\ep_1,\ldots, \ep_n)$ is chosen so that the associated $SO(3)$ bundle has second Stiefel-Whitney class $w_2 \neq 0$ nontrivial, then every fixed point of $\ep \si$ in $\Hom(F_n, SU(2))$ is necessarily irreducible.

 We  view conjugacy classes of fixed points of $\ep \si$ as points in the intersection  $\wh \De \cap \wh \Ga_{\ep \si},$ where
$\wh \De =\De /\conj $ is the quotient of the diagonal $\De \subset C_i^n \times C_i^n$,   and where $\wh \Ga_{\ep \si} = \Ga_{\ep \si} /\conj $ is the quotient of the graph $\Ga_{\ep \si}$ of $\ep \si \colon C_i^n \to C_i^n$.   

If the link $L$ is the closure of a 2-strand braid, as it is for the Hopf link, then $\ep=(-1,-1)$ is the only choice whose associated $SO(3)$ bundle has $w_2 \neq 0$. 
Furthermore, in this case 
the intersection $\wh \De \cap \wh \Ga_{\ep \si}$ takes place in the pillowcase  
$f^{-1}(1)/\conj$, which is defined as the quotient 
\begin{equation}\label{PC-sleepy}
P=\{ (a,b,c,d) \in C_i^4 \mid ab=cd\}/\conj.
\end{equation}

It is well known that $P$  is homeomorphic to $S^2$.  To see this, first conjugate so that $a=i$, then conjugate by elements of the form $e^{i \th}$ to arrange that $b$ lies in the $(i,j)$-circle. A straightforward calculation using the equation $ab=cd$ shows    that $d$ must also lie on the $(i,j)$-circle. Clearly   $c$ is determined by $a,b,d$.    We thus obtain an embedded 2-torus of elements of $C_i ^4$  satisfying $ab=cd$, parameterized by   
$$g(\th_1,\th_2) = (i, e^{k\th_1}i, e^{k(\th_2-\th_1)}i, e^{k\th_2}i)$$
for $\th_1,\th_2 \in [0,2\pi), $ which maps onto $P$.  It is easy to verify that this is a two-to-one 
submersion,   except when  $\theta_1, \theta_2\in \{ 0, \pi \}$.  This realizes $P$ as a quotient of the torus by the hyperelliptic involution.  In particular, this involution is orientation preserving, and away from the four singular points of $P$, we can lift all orientation questions  up to the torus.

 Let $L$ be the  right-handed  Hopf link, which we view as the closure of the braid $\si = \si_1^2 \in B_2$, and suppose $\ep = (-1,-1).$   
The intersection point $\wh \De \cap \wh \Ga_{\ep \si}$ in question is given by the conjugacy class of $g(\pi/2,\pi/2)$, namely the point $[(i,j,i,j)] \in P$.  Thus, in order to pin down the sign of  the Casson-Lin invariant $h_2(L)$, we must determine the orientations of $\wh \De, \; \wh \Ga_{\ep \si},$ and $P$ at this point.

Notice that
\begin{eqnarray*}
\frac{\partial}{\partial \th_1} g(\th_1,\th_2) &=& (0, e^{k\th_1}j, -e^{k(\th_2-\th_1)}j, 0)\\
\frac{\partial}{\partial \th_2} g(\th_1,\th_2) &=& (0, 0, e^{k(\th_2-\th_1)}j, e^{k\th_2}j).
\end{eqnarray*}
Evaluating at $\th_1 = \th_2 = \pi/2$ gives two tangent vectors 
$u_1 :=(0,-i,-j,0)$ and $u_2:=(0,0,j,-i)$  to $C_i ^4$   which span   a complementary subspace in  $\ker df$ to the orbit tangent space.  Therefore, an ordering of these vectors determines an orientation on $P=f^{-1}(1)/\conj$.   

The orbit tangent    space is spanned by  the three tangent vectors 
\begin{eqnarray*}
v_1 :=\left. \frac{\partial}{\partial t}\right|_{t=0} e^{it} (i,j,i,j) &=& (0, k,0,k), \\
v_2 := \left. \frac{\partial}{\partial t}\right|_{t=0} e^{jt} (i,j,i,j) &=& (-k,0,-k, 0),  \\
v_3:= \left. \frac{\partial}{\partial t}\right|_{t=0} e^{kt} (i,j,i,j) &=& (j, -i,j,-i).  \end{eqnarray*}

Then the five vectors $\{u_1,u_2, v_1, v_2, v_3\}$ form a basis for ${\rm ker}\left(df|_{(i,j,i,j)}\right) =T_{(i,j,i,j)} f^{-1}(1).$ We choose vectors $w_1=(k,0,0,0), w_2 = (0,k,0,0), w_3=(j,0,0,0)$ to extend this to a basis  $\{u_1,u_2, v_1, v_2, v_3, w_1, w_2, w_3\}$ for  $T_{(i,j,i,j)}C_i^4.$

 The orientation conventions in the definition of the $h_{N,a}(L)$ \cite{BHa} involve pulling back the orientation from $su(2)=T_1 SU(2)$ by $df$ to obtain a co-orientation for ${\rm ker}\left(df|_{(i,j,i,j)}\right)$. With that in mind,  we compute the action of $df$ 
on $\{w_1, w_2, w_3\}$, namely,  $df(w_1)=-j$, $df(w_2)=i$ and $df(w_3)=k$.  

Notice that the ordered triple $\{df(w_1), df(w_2), df(w_3)\} = \{-j,i,k\}$ gives the same orientation as  the standard basis for $su(2).$  
Thus,    the base-fiber rule gives the co-orientation  $\{w_1, w_2, w_3\}$     on 
$\ker df$,   so    we choose the  orientation $\cO_{\ker df} $ on $\ker df$ such that 
$\cO_{\{w_1, w_2, w_3\}} \oplus \cO_{\ker df} $ agrees with the product orientation on $C_i ^2 \times C_i^2$.

The orientation on the pillowcase $P$ is then obtained by applying the base-fiber rule a second time to the quotient \eqref{PC-sleepy}, using  $\cO_{\ker df}$ to orient $f^{-1}(1)$ and giving the orbit tangent space the orientation induced from that on $SU(2)$ as well. We claim that the basis $\{u_1, u_2\}$ for the tangent space to the pillowcase has the opposite  orientation.   
To see this, we note that $\{v_1, v_2, v_3 \}$ is the fiber orientation for $SO(3)\to f^{-1}(1) \to P$ and compare  
$$S=\{w_1, w_2, w_3, u_1,u_2, v_1, v_2, v_3 \}$$ to the  product orientation on $C_i ^2 \times C_i ^2$.
Using the basis 
 $$ \{(j,0),(k,0),(0,k),(0,i) \}$$  
for $T_{(i,j)} (C_i ^2 )$, we see that    
\begin{eqnarray*}
\be &=& \{(j,0,0,0),(k,0,0,0),(0,k,0,0),(0,i,0,0), \\
&& \quad (0,0,j,0), (0,0,k,0), (0,0,0,k),(0,0,0,i)\}
\end{eqnarray*}
is an oriented basis for    $T_{(i,j,i,j)}C_i^4 = T_{(i,j)}C_i^2 \times T_{(i,j)}C_i^2$ with the product orientation.   

Let $M$ be the matrix   expressing    the vectors in $S$ in terms of the basis $\be$.   Since   
$$M=
\begin{bmatrix}
0&0&1&0&0&0&0&1\\
1&0&0&0&0&0&\matminus1&0\\
0&1&0&0&0&1&0&0\\
0&0&0&\matminus1&0&0&0&\matminus1\\
0&0&0&\matminus1&1&0&0&1\\
0&0&0&0&0&0&\matminus1&0\\
0&0&0&0&0&1&0&0\\
0&0&0&0&\matminus1&0&0&\matminus1
\end{bmatrix},$$
one easily computes that $\det M=-1$,   confirming  our claim that $\{ u_2, u_1\}$ is a positively oriented basis for the pillowcase tangent space.

\bigskip
\begin{figure}[h]
\includegraphics[scale=1.0]{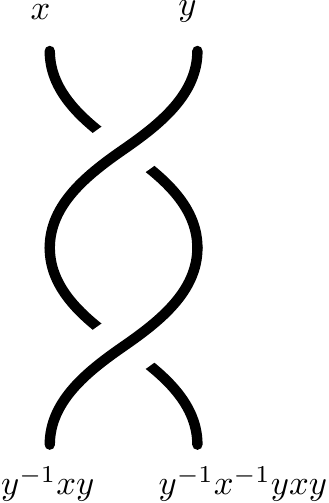} 
\caption{The action of $\si=\si_1^2$ on $F_2=\langle x,y \rangle$}
\label{braid-comp}
\end{figure}

Recall that $L$ is the right-handed Hopf link, which we represent as the closure of the braid $\si = \si_1^2 \in B_2$.
For $\ep=(-1,-1)$,
as in Figure \ref{braid-comp}, one can verify that
$$\ep \si(X,Y)= (-Y^{-1}XY, -Y^{-1}X^{-1}YXY).$$
 Consider the    curve $\al(\th) = (i, e^{k\th} i)$,  passing through the point $(i,j)\in C_i ^2$ when $\th={\pi}/{2}$,  which is transverse to the orbit $[(i,j)]$.  Then $(\al (\th), \al (\th) )$ and $(\al (\th),  \ep \si \circ \al(\th) )$  are curves in 
$\De$ and $\Ga_{\ep \si}$,  respectively, and both are necessarily transverse to the orbit in   $C_i ^4 /\conj$.    Therefore, we can   compare the orientations induced by the parameteriziations    $[(\al (\th), \al(\th) )]$ and  $[(\al (\th), \ep \si \circ \al ( \th) )]$ of   $\wh \De$ and $\wh \Ga_{\ep \si}$  to the pillowcase orientation determined above, namely $\{ u_2, u_1\}$. 
The velocity vectors for  the paths $(\al (\th) , \al (\th) )=(i, e^{k\th}i, i, e^{k\th}i)$ and $(\al(\th) , \ep \si \circ \al (\th) ) = (i, e^{k\th}i, 
-e^{2k\th}i,-e^{3k\th}i )$ at 
$\th={\pi}/{2}$ are given by  $(0,-i,0,-i)=u_1+u_2$  and $(0,-i,2j,-3i)=u_1 + 3u_2$, respectively. 
 
The Casson-Lin invariant is defined as the intersection number $h_2(L)= \langle \wh \De,\; \wh{\Ga}_{\ep \si}\rangle,$ and  in our case
the sign of the unique intersection point in $\wh \De \cap \wh{\Ga}_{\ep \si}$ is determined by comparing the orientation of $\{u_1+u_2, u_1+3u_2\}$ with 
$\{u_2, u_1\}.$  Since  the change of basis matrix $\SmallMatrix{1&\,3\\1&\,1}$
has negative determinant, it follows that  $h_2(L)=-1,$ and this completes the proof of the theorem. \qed

 \bigskip
\noindent
 {\it Acknowledgements.}  The authors would like to thank Eric Harper and Nikolai Saveliev for  many helpful discussions.

\end{document}